\documentclass[12pt]{article}

\usepackage{amssymb,amsmath,epsf}

\newtheorem{lemma}{Lemma}
\newtheorem{proposition}{Proposition}
\newtheorem{theorem}{Theorem}
\newtheorem{remark}{Remark}
\newtheorem{definition}{Definition}

\newenvironment{proof}{\smallskip\noindent{\it Proof.}\hskip \labelsep}
                        {\hfill\penalty10000\raisebox{-.09em}{$\Box$}\par\medskip}
\newfont{\beb}{msbm10 at 12pt}
\newfont{\roc}{eusm10 at 12pt}

\def\R{\mathbb{R}}
\def\C{\mathbb{C}}
\def\Z{\mathbb{Z}}
\def\esf{\mathbb{S}}
\def\T{\mathbb{T}}
\def\M{\mathbb{M}}
\def\l{{\lambda}}
\def\a{{\alpha}}
\def\be{{\beta}}
\def\t{{\theta}}
\def\g{{\gamma}}

\def\rth{{\mathbb{R}^3}}

\begin{document}
\title{The space of doubly periodic minimal tori with parallel ends: Standard examples}
\author{M. Magdalena Rodr\'\i guez\thanks{Research partially supported by a MEC/FEDER grant no. MTM2004-02746.}}
\date{}
\maketitle

\noindent {\sc Abstract.} {\footnotesize We describe a 3-parametric
  family ${\cal K}$ of properly embedded minimal tori with four
  parallel ends in quotients of $\R^3$ by two independent
  translations, which we will call the {\it standard examples.}  These
  surfaces generalize the examples given by Karcher, Meeks and
  Rosenberg in~\cite{ka4,ka6,mr3}. ${\cal K}$ can be endowed with a
  natural structure of a self-conjugated 3-dimensional real analytic
  manifold diffeomorphic to $\R\times\left(\R^2-\{(\pm 1,0)\}\right)$
  whose degenerate limits are the catenoid, the helicoid, the simply
  and doubly periodic Scherk minimal surfaces and the Riemann minimal
  examples. P\'{e}rez, Rodr\'\i guez and Traizet~\cite{PeRoTra1}
  characterize ${\cal K}$ in the following sense: If $M$ is a properly
  embedded minimal torus in a quotient of $\R^3$ by two independent
  translations with any number of parallel ends, then $M$ is a finite
  covering of a standard example.}

\section{Introduction}
\label{secintrod}
Scherk~\cite{sche1} presented in 1935 the first properly embedded minimal surface\footnote{Unless explicitly mentioned, all surfaces in the paper are supposed to be connected and orientable.} in $\R^3$, invariant by two linearly independent translations (we will shorten by saying a {\it doubly periodic minimal surface}).
This surface is known as {\it Scherk's first surface}, and fits naturally into a 1-parameter family ${\cal F}=\{F_\t\}_\t$ of examples, called doubly periodic Scherk minimal surfaces.
In the quotient by its more refined period lattice (i.e. the period lattice generated by its shortest period vectors), each $F_\t$ has genus zero and four asymptotically flat annular ends: two top and two bottom ones, provided that the period lattice is horizontal.
This kind of annular ends are called {\it Scherk-type ends}.
The parameter $\t$ in this family ${\cal F}$ is the angle between top and bottom ends of $F_\t$.
We can clearly consider the quotient of these $F_\t$ by less refined period lattice to have two top and $2k$ bottom ends for any natural $k$, keeping genus zero in the quotient.
Lazard-Holly and Meeks~\cite{ml1} proved that these are the only possible examples in this setting;
i.e., if the quotient of a doubly periodic minimal surface $M\subset\R^3$ has genus zero, then $M$ must be a doubly periodic Scherk minimal surface up to translations, rotations and homotheties.
Moreover, the angle map $\t:{\cal F}\to(0,\pi)$ is a diffeomorphism.
Hence the moduli space of
 properly embedded minimal surfaces with genus zero in $\T\times\R$, $\T$ a flat torus, is diffeomorphic to $(0,\pi)$ after identifying by rotations, translations and homotheties.

In 1988, Karcher~\cite{ka4} defined another $1$-parameter family of doubly periodic minimal surfaces, called {\it toroidal halfplane layers}, with genus one and four Scherk-type parallel ends in its smallest fundamental domain (these examples will be denoted as $M_{\t,0,0}$ in Section~\ref{secstandardexamples}).
Furthermore, he exposed two distinct $1$-parameter deformations of each toroidal halfplane layer, obtaining other doubly periodic minimal tori with parallel ends
(denoted as $M_{\t,\a,0}$ and $M_{\t,0,\be}$, with $\be<\t$, in Section~\ref{secstandardexamples}).
We generalize these Karcher's examples in Section~\ref{secstandardexamples}, by obtaining a $3$-parameter family.

\begin{theorem}
\label{thmK}
There exists a $3$-parameter family ${\cal K}=\{M_{\t,\a,\be}\}_{\t,\a,\be}$ of properly embedded doubly periodic minimal surfaces with genus one and four parallel ends in the quotient by their more refined period lattices.
This family ${\cal K}$ can be endowed with a natural structure of a real analytic $3$-dimensional manifold with the topology of the uniform convergence on compact sets. Furthermore:
\begin{enumerate}
\item The isometry group of any surface $M_{\t,\a,\be}\in {\cal K}$ is isomorphic to $(\Z /2\Z)^2$, $(\Z /2\Z)^3$ or $(\Z /2\Z)^4$, depending on the values of $\a,\be$, and it contains an orientation reversing involution without fixed points, producing a quotient Klein bottle with $2$ parallel ends.
\item ${\cal K}$ is {\it self-conjugate}, in the sense that the conjugate surface\footnote{Two minimal surfaces $M_1,M_2\subset \R^3$ are {\it conjugate} if the coordinate functions of $M_2$ are harmonic conjugate to the coordinate functions of $M_1$.} of any example in ${\cal K}$ also belongs to ${\cal K}$.
\item The possible limits of surfaces in ${\cal K}$ are the catenoid, the
helicoid, any singly or doubly periodic Scherk minimal surface, any Riemann minimal example or another surface in ${\cal K}$.
\end{enumerate}
\end{theorem}

We refer to the examples $M_{\t,\a,\be}$ in ${\cal K}$ as {\it standard examples}.
It is clear that we can consider quotients of the standard examples by less refined period lattices to have $4k$ ends for any natural $k$, keeping genus one.
P\'erez, Rodr\'\i guez and Traizet~\cite{PeRoTra1} have proved that these are the only possible examples in this setting.

\begin{theorem}\cite{PeRoTra1}
\label{thmuniqueness}
If $M\subset\R^3$ is a doubly periodic minimal surface with parallel ends and genus one in the quotient, then $M$ must be a standard example in ${\cal K}$ up to translations, rotations and homotheties.
\end{theorem}

Meeks and Rosenberg~\cite{mr3} developed a general theory for doubly periodic minimal
surfaces having finite topology in the quotient, and used an approach of minimax type to find theoretically some new examples with parallel ends and genus one in the quotient, besides those given by Karcher.
After studying in detail the surfaces in ${\cal K}$, the uniqueness Theorem~\ref{thmuniqueness} assures that Meeks and Rosenberg's examples are nothing but $M_{\t,0,\be}$, for $\be<\t$.
Thus at least two of the most symmetric 1-parameter families in ${\cal K}$ were known by Karcher~\cite{ka4,ka6} and by Meeks and Rosenberg~\cite{mr3} (although our approach here is different from theirs). For this reason, the surfaces in ${\cal K}$ also appear sometimes in the literature as {\it KMR examples}.

We will construct all standard examples as branched coverings of the sphere $\esf^2$ by their Gauss maps. The {\it spherical configuration} of a standard example, defined as the position in $\esf^2$ of the branch values of its Gauss map, allows us to read all the information concerning the minimal surface, see Section~\ref{secstandardexamples}.
Besides giving a unified method to produce all standard examples and studying their geometry, our motivation for writing this paper was to study the topology of ${\cal K}$.

\begin{theorem}
\label{thmtopology}
The space ${\cal K}$ of properly embedded minimal surfaces with genus one and parallel ends in $\T\times\R$, $\T$ a $2$-dimensional flat torus, is diffeomorphic to $\R\times(\R^2-\{(\pm 1,0)\})$.
\end{theorem}

The proof of Theorem~\ref{thmtopology} is inspired by the arguments of P\'erez, Traizet and the author~\cite{PeRoTra1} to prove Theorem~\ref{thmuniqueness} (they follow the ideas of Meeks, P\'erez and Ros~\cite{mpr1}).
We model the family ${\cal K}$ of standard examples as an analytic subset in a complex manifold ${\cal W}$ of finite dimension (roughly, ${\cal W}$ consists of all admissible Weierstrass data for our problem).
In the boundary of ${\cal K}$ in ${\cal W}$, we can find the $1$-parameter family ${\cal S}$ of singly periodic Scherk minimal surfaces~\cite{ka4,sche1}.
We consider the classifying map $C:\widetilde{\cal K}\to\Lambda=\R^+\times\esf^1\times\R$, defined on $\widetilde{\cal K}={\cal K}\cup{\cal S}$,€€ which associates to each surface in $\widetilde{\cal K}$ two geometric invariants: its period at the ends and its flux along a nontrivial homology class with vanishing period vector.
Theorem~\ref{thmtopology} is a simple consequence of the following statements:
\begin{enumerate}
    \item $C:\widetilde{\cal K}\to\Lambda$ is a proper map.
    \item $C:\widetilde{\cal K}\to\Lambda$ is a local diffeomorphism.
    \item There exists $x\in\Lambda$ such that $C^{-1}(x)$ consists of only one surface in $\widetilde{\cal K}$.
    \item $C({\cal S})$ is a proper, divergent curve in $\Lambda$.
\end{enumerate}
\noindent

The paper is organized as follows.
In Section~\ref{secstandardexamples} we study the family ${\cal K}$ of standard examples.
Section~\ref{secclassifyingmap} is devoted to introduce the space ${\cal W}$ of admissible Weierstrass data and the classifying map $C$ that we use as a tool to demonstrate Theorem~\ref{thmtopology}, and we prove that $C$ is a proper map.
The goal of Section~\ref{secDiff} is to prove the second statement above; i.e. $C$ is a local diffeomorphism.
Finally, it can be found in Section~\ref{secKnotsimplyconnected} the proof of Theorem~\ref{thmtopology}.

I sincerely want to thank Joaqu\'\i n P\'erez for his invaluable hepl along these years, and for leading me through this work.

\section{Standard examples  (proof of Theorem~\ref{thmK})}
\label{secstandardexamples}
We dedicate this section to introduce the $3$-parameter family ${\cal K}$ of standard examples appearing in Theorem~\ref{thmK}, to which the uniqueness Theorem~\ref{thmuniqueness} applies.
First, let us point out some general facts.
Let $\widetilde{M}\subset \R^3$ be a doubly periodic minimal surface with period lattice~${\cal P}$.
Such $\widetilde{M}$ induces a properly embedded minimal surface $M=\widetilde{M}/{\cal P}$ in the complete flat $3$-manifold $\R^3/{\cal P}= \T\times\R$, where $\T$ is a $2$-dimensional flat torus.
Reciprocally, if $M\subset \T\times \R$ is a properly embedded nonflat minimal surface, then its lift $\widetilde{M}\subset\R^3$ is a connected doubly periodic minimal surface, by the Strong Halfspace Theorem of Hoffman and Meeks~\cite{hm10}.
Assume that the topology of $M$ is a finitely punctured torus and that its ends are parallel.
Then Meeks and Rosenberg~\cite{mr3} ensure that $M$ has finite total curvature and $4k$ Scherk-type ends, for some natural $k$.
Therefore $M$ is conformally equivalent to a torus $\M$ minus $4k$ punctures.
If we consider ${\cal P}$ to be the more refined period lattice of $\widetilde{M}$, then Theorem~\ref{thmuniqueness} implies that $k=1$.

Since $M$ has finite total curvature, its Gauss map $g$ extends meromorphically to $\M$.
After a rotation so that the ends of $M$ are horizontal, $g$ takes values $0,\infty$ at the punctures, and the third coordinate function $h$ (which is not well-defined on $M$) defines an univalent holomorphic $1$-form $dh$ on $\M$, which we will call the height differential.
Meeks and Rosenberg~\cite{mr3} proved that one of the meromorphic differentials $g\,dh,\frac{dh}{g}$ has a simple pole at each puncture.
As $dh$ has no zeros on $\M$ (it has no poles), we conclude that $g$ is unbranched at the ends, and has degree two.
The Riemann-Hurwitz formula implies that the total branching number of $g$ is four.

Any standard example will be given in terms of the branch values of its Gauss map, which will consist of two pairs of antipodal points $D,D',D'',D'''$ in the sphere $\esf^2$.
We label those points so that $D''=-D$, $D'''=-D'$.
Since the Gauss map is unbranched at the ends (which are horizontal), we impose these branch values to be different from the North and South Poles.
We denote by $e\subset\esf^2$ the equator that contains $D,D',D'',D'''$ and by $P\in e$ the point that bisects the angle $2\t$ between $D$ and $D'$, $\t \in(0,\frac{\pi}{2})$.
We will call a {\it spherical configuration} to any set $\{D,D',D'',D'''\}$ as above.

\subsection{Toroidal halfplane layers $M_{\t,0,0}$}
\label{subsecMt}
\begin{figure}
\begin{center}
\epsfysize=5cm 
\epsffile{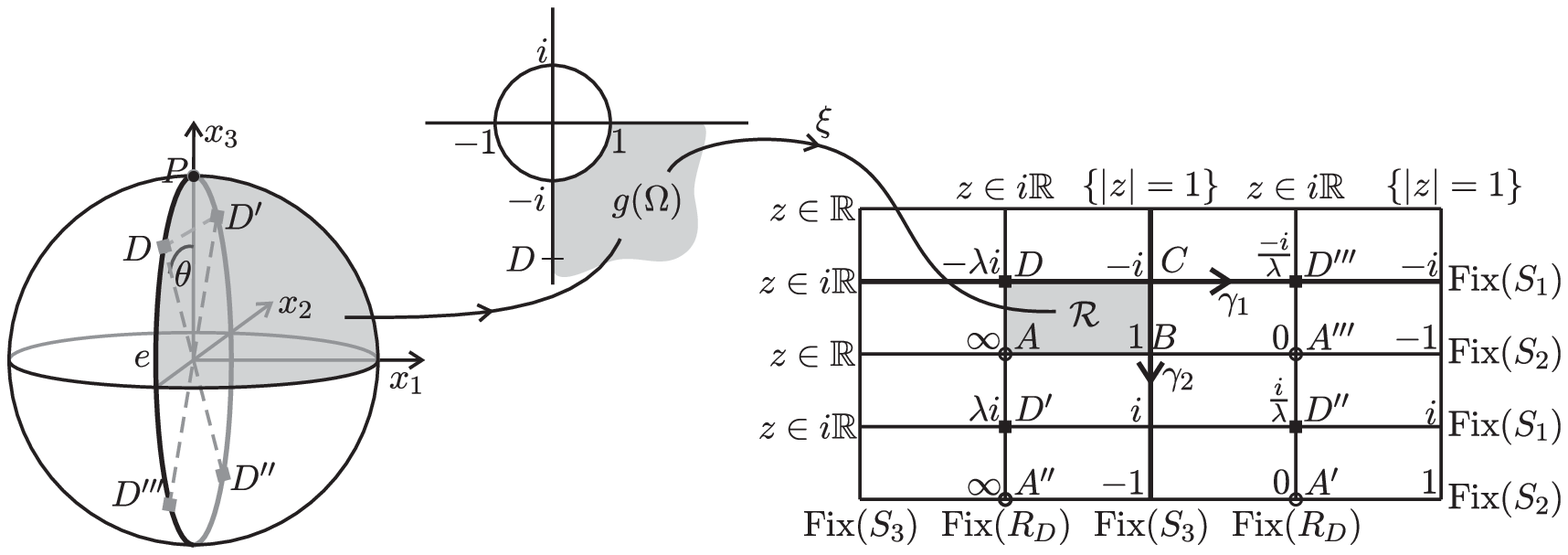}
\end{center}
\caption{Left: Spherical configuration of $M_{\t,0,0}$. Center:
The biholomorphism $\xi$ between the shaded regions. Right: The
conformal torus $\Sigma_\t$, where Fix$(\bullet)$ denotes fixed
point set.}
\label{M_t00}
\end{figure}

With the notation above, given $\t \in(0,\frac{\pi}{2})$ we set the equator $e$ to be the inverse image of the imaginary axis $i\overline\R\subset\overline\C $ through the stereographic projection from the North Pole, and ${P=(0,0,1)}$, see Figure~\ref{M_t00} left.
After stereographic projection we have $D=-\l i$, with $\l =\l(\t)=\cot\frac{\t}{2}$, and the remaining branch values of the Gauss map of the example $M_{\t,0,0}$ we are constructing are the four roots of the polynomial $(z^2+\l^2)(z^2+\l^{-2})$.
Thus the underlying conformal compactification of the potential surface $M_{\t,0,0}$ is the rectangular torus
\[
\Sigma_\t=\left\{(z,w)\in \overline{\C}^2\ |\ w^2=(z^2+\l^2)(z^2+\l^{-2})\right\}.
\]
The degree two extended Gauss map of $M_{\t,0,0}$ is $g(z,w)=z$, the punctures correspond to $(0,\pm 1),(\infty,\pm \infty)\in\Sigma_\t$, and the height differential must be $dh=\mu\frac{dz}{w}$ for certain ${\mu=\mu(\t)\in\C^*}$.

We consider $\mu \in \R^*$. Then the set $\{(z,w)\ |\ |z|=1\}$ corresponds on $M_{\t,0,0}$ to two closed horizontal geodesics which are the fixed point set of reflection symmetries $S_3$ in two horizontal planes (the reflection in both planes induce the same isometry $S_3$ of the quotient surface);
the set $\{(z,w)\ |\ z\in \R \}$ corresponds on $M_{\t,0,0}$ to four geodesics traveling from a zero to a pole of the Gauss map $g$, which are the fixed point set of a reflection symmetry $S_2$ across two planes orthogonal to the $x_2$-axis;
and the set $\{(i t,w)\ |\ t\in\R, \l^{-1}\leq|t|\leq\l \}$ corresponds to two geodesics which are the fixed point set of a reflection symmetry $S_1$ in a vertical plane orthogonal to the $x_1$-axis.
The later geodesics cut orthogonally four straight lines parallel to the $x_1$-axis and contained in $M_{\t,0,0}$, which correspond to the set $\{(i t,w)\ |\ |t|\leq\l^{-1}\ \mbox{or }|t|\geq\l \}$.
We will denote by $R_D$ the $\pi$-rotation around any such straight line, see Figure~\ref{StandardExamples} left.

We now construct a different model of $\Sigma_\t$, as a quotient of the $\xi$-plane $\C$ over a rectangular lattice.
Let $\Omega \subset \Sigma_\t$ be one of the two components of $g^{-1}(\{|z|>1,-\frac{\pi}{2}<\arg(z)<0\})$. $\Omega$ is topologically a disk and its boundary contains the branch point corresponding to the branch value $D$ of $g$ and one of the ends corresponding to a pole of $g$.
Let ${\cal R}$  be an open rectangle in the $\xi$-plane of consecutive vertices $A,B,C,D\in \C $ with the segment $\overline{AB}$ being horizontal, such that there exists a biholomorphism $\xi:\{|z|>1,-\frac{\pi}{2}<\arg(z)<0\}\to {\cal R}$ with boundary values $\xi(\infty)=A$, $\xi(1)=B$, $\xi(-i)=C$ and $\xi(-\l i)=D$.
Then the composition of $g$ with $\xi$ defines a biholomorphism between $\Omega$ and ${\cal R}$.
After symmetric extension of this biholomorphism across the boundary curves we will get a biholomorphism from $\Sigma_\t$ to the quotient of the $\xi$-plane modulo the translations given by four times the sides of the rectangle ${\cal R}$.
We abuse the notation by labeling also $D,D',D'',D'''$ the points of the $\xi$-plane that correspond to the branch values of $g$.
The deck transformation $(z,w)\stackrel{\cal D}{\mapsto}(z,-w)$ of $\Sigma_\t$ corresponds in the $\xi$-plane to the $\pi$-rotation about the branch points of $g$.
It will be also useful to see $\Sigma_\t$ as a branched $2:1$ covering of $\overline\C$ through the map $(z,w)\mapsto z$, i.e. two copies $\overline{\C}_1,\overline\C_2$ of $\overline\C$ glued along common cuts from $D$ to $D'$ and from $D''$ to $D'''$, both contained in the imaginary axis.

In the $\xi$-plane model of $\Sigma_\t$, $S_3$ corresponds to the reflection across the line passing through $B,C$
(or across its parallel line after translation by half a horizontal period, see Figure~\ref{M_t00} right);
$S_2$ is the reflection across the line passing through $A,B$
(or across its parallel line after translation by half a vertical period);
$S_1$ is the reflection across the line passing through $D,D'''$ (or through $D',D''$);
and $R_D$ is the reflection across the line passing through $D,D'$ (or through $D'',D'''$).
It is easy to see that Iso$(M_{\t,0,0})$ coincides with the group of conformal transformations of the underlying conformal torus $\Sigma_\t$, which is isomorphic to $(\Z /2\Z)^4$ with generators $S_1,S_2,S_3,R_D$.

\begin{remark}
\label{notaisometria}
All surfaces $M_{\t,\a,\be}\in {\cal K}$ to be defined will have the same conformal compactification $\Sigma_\t$ as $M_{\t,0,0}$.
So from now on the sixteen elements in Iso$(M_{\t,0,0})$ generated by $S_1,S_2,S_3,R_D$ will be seen as conformal transformations of $\Sigma_\t$.
Those that leave invariant the distribution of zeros and poles of the Gauss map of $M_{\t,\a,\be}$ will be precisely the isometries of this last surface.
\end{remark}

Concerning the period problem for $M_{\t,0,0}$, let $\g_1,\g_2$ be the simple closed curves in $\Sigma_\t$ obtained as quotients of the horizontal and vertical lines in the $\xi$-plane passing through $D,D'''$ and through $C,B$ respectively (see Figure~\ref{M_t00} right).
Clearly $\{\g_1,\g_2\}$ is a basis of $H_1(\Sigma_{\t},\Z)$.
We normalize so that $\int_{\g_2}dh=2\pi i$, which determines $dh$ or, equivalently, the value of $\mu$,
\begin{equation}
\label{eq:mu}
\mu=\mu(\t) = \frac{\pi\csc\t}{\mbox{\roc K}(\sin^2\t)},
\end{equation}
where $\mbox{\roc K}(m)=\int_0^{\frac{\pi}{2}}\frac{1}{\sqrt{1-m\sin^2u}}\, du$ , $0<m<1$, is the complete elliptic integral of the first kind.
With this choice of $\mu$, the period and flux vectors of $M_{\t,0,0}$ along a small loop $\g_A$ around the end $A=(\infty,+\infty)\in \Sigma_\t$ of $M_{\t,0,0}$ are respectively
\begin{equation}
\label{eq:periodendsmuysim}
P_{\g_A}=(0,\pi \mu,0)\quad{\rm and}\quad F_{\g_A}=(\pi \mu,0,0).
\end{equation}
The remaining ends of $M_{\t,0,0}$ are
\[
A'=(S_1\circ S_2\circ S_3)(A),\quad A''={\cal D}(A)=(S_1\circ R_D)(A),\quad A'''={\cal D}(A'),
\]
see Figure~\ref{M_t00} right.
From the behavior of the Weierstrass form $\Phi=\left(\frac{1}{2}(\frac{1}{g}-g),\frac{i}{2}(\frac{1}{g}+g),1\right)dh$ under pullback by
$S_1,S_2,S_3,R_D$, one obtains that
\begin{equation}
\label{eq:residuesallends}
\mbox{Res}_A\Phi=-\overline{\mbox{Res}_{A'}\Phi}=-\mbox{Res}_{A''}\Phi=\overline{\mbox{Res}_{A'''}\Phi},
\end{equation}
where $\mbox{Res}_X$ denotes the residue at the point $X\in\Sigma_\t$.
Note that (\ref{eq:periodendsmuysim}) and (\ref{eq:residuesallends}) determine completely the periods and fluxes at $A',A'',A'''$:
\begin{equation}
\label{eq:finales}
P_{\g_A}=P_{\g_{A'}}=-P_{\g_{A''}}=-P_{\g_{A'''}}\qquad\mbox{ and }\qquad F_{\g_A}=-F_{\g_{A'}}=-F_{\g_{A''}}=F_{\g_{A'''}} .
\end{equation}

Similar arguments imply that the periods and fluxes along the homology basis are
\begin{equation}
\label{eq:periodhomologymuysim}
\begin{array}{lll}
P_{\g_1}=( 0,0,f_1),&\qquad &F_{\g_1}=-F_{\g_A}=(-\pi \mu,0,0),\\
P_{\g_2}=(0,0,0),     &\qquad &F_{\g_2}=(0,0,2\pi),
\end{array}
\end{equation}
where
\begin{equation}
\label{eq:f1}
f_1=f_1(\t)=-4\mu\int_1^{\l}\frac{dt}{\sqrt{(t^2-\l^{-2})(\l^2-t^2)}}<0.
\end{equation}
From equations (\ref{eq:periodendsmuysim}), (\ref{eq:finales}) and (\ref{eq:periodhomologymuysim})
we conclude that $M_{\t,0,0}$ is a complete immersed minimal surface invariant by the rank two lattice generated by $P_{\g_A},P_{\g_1}$.
Moreover, $M_{\t,0,0}$ has genus one and four horizontal Scherk-type ends in the quotient, and can be decomposed in 16 congruent disjoint pieces.
Karcher~\cite{ka4} proved that each of these pieces is the conjugate surface of certain Jenkins-Serrin graph defined on a convex domain.
In particular, $M_{\t,0,0}$ is embedded.

\begin{figure}
\begin{center}
\epsfysize=35mm 
\epsffile{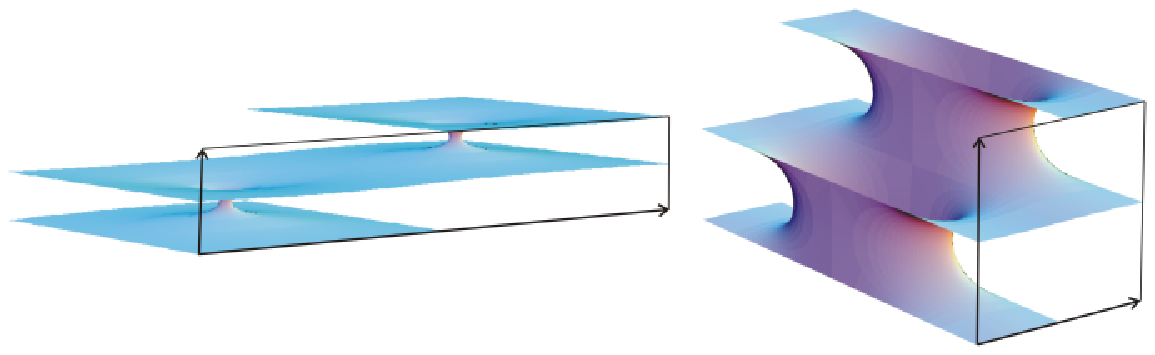}
\end{center}
\caption{$M_{\t,0,0}$ for $\t =\frac{\pi}{50}$ (left) and for $\t =\frac{24 \pi}{50}$ (right).}
\label{figureLimits}
\end{figure}

Next we study the limit surfaces of the examples in the family $\{M_{\t,0,0}\ |\ \t\in (0,\frac{\pi}{2})\}$.
When $\t$ goes to zero, the function $\l(\t)$ diverges to $+\infty$.
After changing variables ${(z,w)\in \Sigma_\t}$ for $(z,w_1)$ with $w_1\l(\t)=w$, it is easy to see that $\Sigma_\t$ degenerates as $\t\to 0^+$ into two spheres $\{(z,w_1)\ |\ w_1^2=z^2\}$.
The limiting Gauss map of $M_{\t,0,0}$ as $\t\to 0^+$ is $g(z,w_1)=z$
and the height differential $dh$ of $M_{\t,0,0}$ converges smoothly to $\frac{dz}{w_1}=\pm \frac{dz}{z}$.
Hence, when ${\t\to 0^+}$, the example $M_{\t,0,0}$ converges smoothly to two vertical catenoids with flux $(0,0,2\pi)$, see Figure~\ref{figureLimits} left.

If $\t\to\frac{\pi}{2}^-$, then $\l(\t)\to 1$ and $\Sigma_\t$ degenerates into two spheres $\{(z,w)\ |
\ w^2=(z^2+1)^2\}$. In this case, the limiting Gauss map is $g(z,w)=z$ and the height differential collapses to zero because the limit of $\mu(\t)$ when $\t \to \frac{\pi}{2}^-$ vanishes.
After scaling, it holds that $\frac{1}{\mu(\t)}dh\to \pm \frac{dz}{z^2+1}$ as $\t \to \frac{\pi}{2}^-$.
Therefore, after blowing up, $M_{\t,0,0}$ converges smoothly as $\t \to \frac{\pi}{2}^-$
to two doubly periodic Scherk minimal surfaces with two horizontal and two vertical ends, see Figure~\ref{figureLimits} right.

\subsection{The examples $M_{\t,\a,\be}$}
\label{subsecMtab}
\begin{figure}
\begin{center}
\epsfysize=43mm 
\epsffile{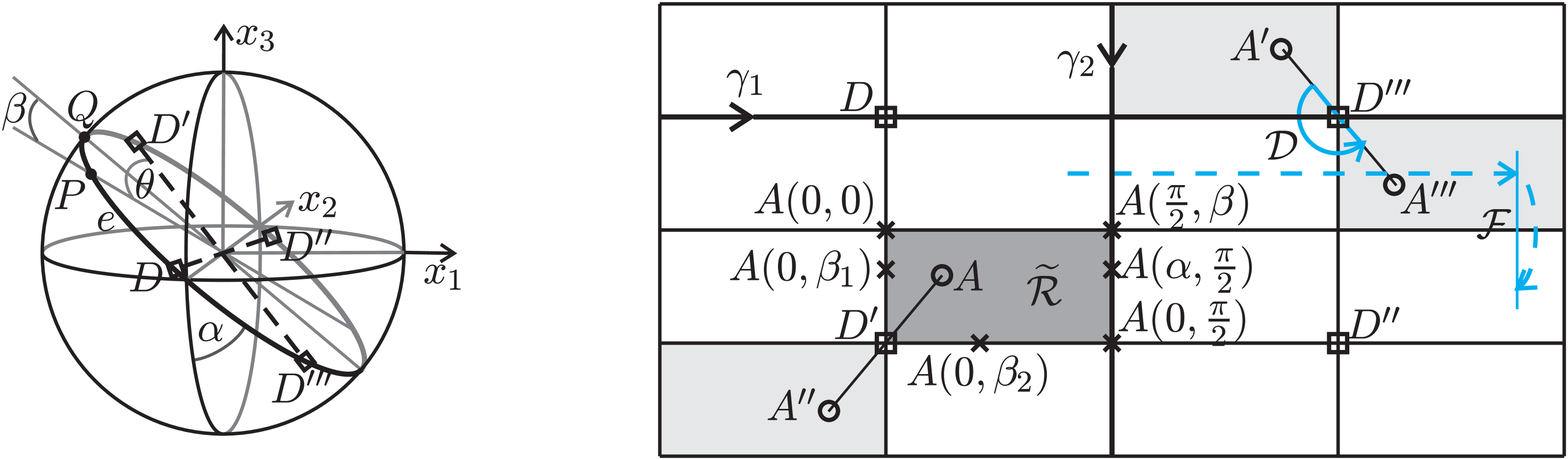}
\end{center}
\caption{Left: Spherical configuration of $M_{\t,\a,\be}$. Right: The $\xi$-plane model of $\Sigma_\t$, where the dotted line represents the isometry ${\cal F}$, $\a\in(0,\frac{\pi}{2})$, $\be\in[0,\frac{\pi}{2}]$ and $0<\be_1<\t<\be_2<\frac{\pi}{2}$.}
\label{ConfEsf}
\end{figure}

Given $\t \in(0,\frac{\pi}{2})$, $\a \in[0,\frac{\pi}{2}]$ and $\be \in [0,\frac{\pi}{2}]$ with $(\a,\be)\neq(0,\t)$, we
consider the equator $e$ to be the rotated image of the imaginary axis in the sphere by angle $\a$ around the $x_2$-axis. If we denote by $Q$ the rotated point by angle $\a$ around the $x_2$-axis of the North Pole, then our new point $P$ will be the rotation of $Q$ by angle $\be$ along $e$, see Figure~\ref{ConfEsf} left.
Note that when $(\a,\be)=(0,\t)$, then $D'$ coincides with the North Pole, which is not allowed in this setting.
Also note that the spherical configuration $\{D,D',D'',D'''\}$ associated to $\t,\a,\be$ is nothing but the rotated image of that of $M_{\t,0,0}$ by the M\"{o}bius transformation $\phi$ corresponding to the composition of the rotation of angle $\be$ around the $x_1$-axis with the rotation of angle $\a$ around the $x_2$-axis.
Consequently, we define the Gauss map $g=g_{\t,\a,\be}$ of the standard example $M_{\t,\a,\be}$ we want to construct as $g= \phi\circ g_{\t,0,0}$, i.e.
\[
g(z,w)=\frac{\sigma\,z +\delta}{i(\overline{\sigma}-\overline{\delta} z)},
\]
for $(z,w)\in \Sigma_\t$, where $\sigma=\cos(\frac{\a +\be}{2})+i \cos(\frac{\a -\be}{2})$ and
$\delta=\sin(\frac{\a -\be}{2})+i \sin(\frac{\a +\be}{2})$.
Since $g$ depends analytically of $\a,\be$, the same holds for its zeros and poles.
We will denote by  $\{A,A',A'',A'''\}=g^{-1}(\{0,\infty\})$ the ends of $M_{\t,\a,\be}$, understanding that each zero or pole of $g$ is defined by analytical continuation of the corresponding zero or pole of $g_{\t,0,0}$.
Choosing the same homology class $[\g_2]\in H_1(\Sigma_\t,\Z)$ as in Subsection~\ref{subsecMt}, we obtain that the height differential of $M_{\t,\a,\be}$ is $dh=\mu\frac{dz}{w}$, with $\mu=\mu(\t)$ as in (\ref{eq:mu}). Thus, the Weierstrass data of $M_{\t,\a,\be}$ coincides with those of $M_{\t,0,0}$ when $\a=\be =0$.

The group Iso$(M_{\t,\a,\be})$ of isometries of the induced metric by $(g,dh)$ always contains the deck transformation ${\cal D}=S_1\circ R_D$ (we follow the notation in Subsection~\ref{subsecMt}, see Remark~\ref{notaisometria}).
Furthermore, the antipodal map in $\esf^2$ leaves invariant the spherical configuration of $M_{\t,\a,\be}$, so Iso$(M_{\t,\a,\be})$ also contains two antiholomorphic involutions without fixed points, ${\cal E}$ and ${\cal F}={\cal E}\circ{\cal D}$.
It is straightforward to check that we can label ${\cal E}=S_1\circ S_2\circ S_3$, and hence ${\cal F}=R_D\circ S_2\circ S_3$.
This information is enough to solve the period problem for $M_{\t,\a,\be}$.

The period and flux vectors of $M_{\t,\a,\be}$ at the end $A$ are given by
\begin{equation}
\label{eq:periodendsnosim}
P_{\g_A}=\pi \mu \sin\t \left( i\, E(\t,\a,\be),0\right),\quad
F_{\g_A}=\pi \mu \sin\t \left( E(\t,\a,\be),0\right),
\end{equation}
where we have used the identification of $\R^3$ with $\C \times \R $ by $(a,b,c)\equiv(a+ib,c)$, and
\[
E(\t,\a,\be)=\frac{1}{\sqrt{\sin^2\t \cos^2\a+(\sin\a \cos\be -i\sin \be)^2}}.
\]
The periods and fluxes at the remaining ends $A'={\cal E}(A)$, $A''={\cal D}(A)$ and $A'''={\cal F}(A)$ can be obtained from the equations in (\ref{eq:finales}), which are still valid.

We choose the homology classes $[\g_1],[\g_2]\in H_1(\Sigma_\t,\Z)$ as in Subsection~\ref{subsecMt}
(note that we can even take the same curve representatives $\g_1,\g_2$ as in the case $\a=\be=0$ except when $\a=\frac{\pi}{2}$ or $\be=\frac{\pi}{2}$).
In particular, the third coordinate $(P_{\g_1})_3$ of the period of $M_{\t,\a,0}$ along $\g_1$ equals $f_1$, given by (\ref{eq:f1}),
so $P_{\g_A},P_{\g_1}$ are linearly independent.
It also holds
\begin{equation}
{\cal E}^*\Phi=-\overline{\Phi},\qquad {\cal E}_*\g_1=-\g_1-\g_A-\g_{A'''},\qquad {\cal E}_*\g_2=\g_2,
\label{eq:calEenhomologia}
\end{equation}
where $\Phi$ denotes the Weierstrass form for $M_{\t,\a,\be}$.
Equalities in (\ref{eq:calEenhomologia}) and (\ref{eq:residuesallends}) imply
$\overline{\int_{\g_1}\Phi}=\int_{\g_1}\Phi +\int_{\g_A}\Phi -\overline{\int_{\g_A}\Phi}$
and  $\int_{\g_2}\Phi= -\overline{\int_{\g_2}\Phi}$, from which we deduce
\begin{equation}
\label{eq:Perg2Mtab}
F_{\g_1}=-F_{\g_A}\qquad \mbox{and} \qquad P_{\g_2}=(0,0,0).
\end{equation}
All of these facts imply that $M_{\t,\a,\be}$ is a complete immersed minimal torus invariant by the rank two lattice generated by $P_{\g_A},P_{\g_1}$, which has four horizontal Scherk-type ends in the quotient.
Since $M_{\t,0,0}$ is embedded and the heights of the ends of $M_{\t,\a,\be}$ depend continuously on $(\a,\be)$, which are in the connected set $[0,\frac{\pi}{2}]^2-\{(0,\t)\}$, we deduce that $M_{\t,\a,\be}$ is embedded outside a fixed compact set.
This fact together with a standard application of the Maximum Principle ensure that $M_{\t,\a,\be}$ is embedded for all values of $\t,\a,\be$.

We next discuss what is the list of isometries of $M_{\t,\a,\be}$ for different values of $\t,\a,\be$.
As we mentioned above, Iso$(M_{\t,\a,\be})$ always contains the subgroup $\{\mbox{identity},{\cal D},{\cal E},{\cal F}\}$, which is isomorphic to $(\Z /2\Z)^2$ with generators ${\cal D}$, ${\cal F}$.
The deck transformation ${\cal D}$ represents in $\R^3$ a central symmetry about any of the four branch points of $g$, and ${\cal F}$ consists of a translation by $\frac{1}{2}(P_{\g_A}+P_{\g_1})$.
In particular, the ends of $M_{\t,\a,\be}$ are equally spaced.
If $0<\be<\frac{\pi}{2}$ and $0<\a<\frac{\pi}{2}$, then the puncture $A=A(\a,\be)$ lies on the open rectangle $\widetilde{\cal R}=S_2({\cal R})$, see Figure~\ref{ConfEsf} right.
By Remark~\ref{notaisometria}, Iso$(M_{\t,\a,\be})$ does not contain either
$S_1, S_2, S_3, R_D$ for these values of $\t,\a,\be$,
and so Iso$(M_{\t,\a,\be})=\{\mbox{identity},{\cal D},{\cal E},{\cal F}\}$.
Thus it remains to study the special cases $\a\in\{0,\frac{\pi}{2}\}$ and $\be\in\{0,\frac{\pi}{2}\}$.

\begin{figure}
\begin{center}
\epsfysize=38mm 
\epsffile{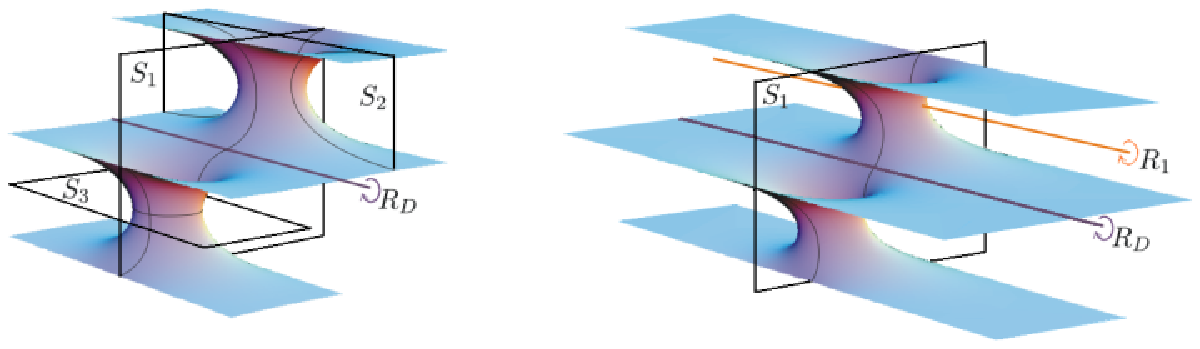}
\end{center}
\caption{Left: $M_{\t,0,0}$ for $\t=\frac{\pi}{4}$. 
Right: $M_{\t,0,\be}$ for $\t=\frac{\pi}{4}$ and $\be=\frac{\pi}{8}$.}
\label{StandardExamples}
\end{figure}

\begin{figure}
\begin{center}
\epsfysize=45mm 
\epsffile{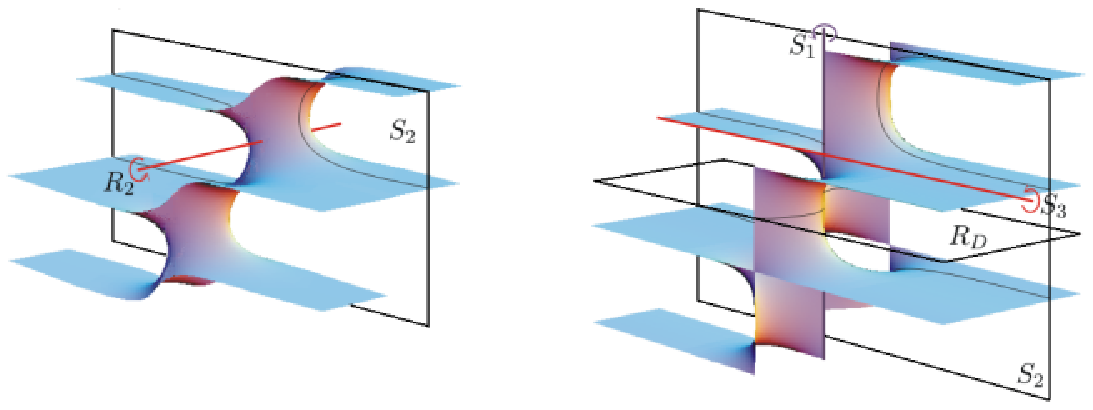}
\end{center}
\caption{Left: $M_{\t,\a,0}$ for $\t=\a=\frac{\pi}{4}$. 
Right: $M_{\t,\a,0}$ for $\t=\frac{\pi}{4}$ and $\a=\frac{\pi}{2}$.}
\label{StandardExamples2}
\end{figure}

\begin{enumerate}
\item The case $\a=\be=0$ was studied in Subsection~\ref{subsecMt}.

\item Suppose that $\a=0$ and $0<\be<\frac{\pi}{2}$, $\be\neq\t$.
In the $\xi$-plane model of $\Sigma_\t$, the puncture~$A$ moves vertically from its original position at the upper left corner of $\widetilde{\cal R}$ when $\be=0$ downwards until collapsing for $\be =\t$ with the branch point $D'$; next it goes on moving horizontally to the right until reaching the lower right corner of $\widetilde{\cal R}$ for $\be=\frac{\pi}{2}$, see Figure~\ref{ConfEsf} right.
The group of isometries Iso$(M_{\t,0,\be})$ is isomorphic to $(\Z /2\Z)^3$ with generators $S_1,R_D,R_1=S_2\circ S_3$.
Here $S_1$ represents in $\R^3$ (as in Subsection~\ref{subsecMt}) a reflection symmetry across a plane orthogonal to the $x_1$-axis, and $R_1$ corresponds to a $\pi$-rotation in $\R^3$ around a line parallel to the $x_1$-axis that cuts the surface orthogonally.
When $0<\be<\t$ (resp. $\t<\be<\frac{\pi}{2}$), $M_{\t,0,\be}$ contains four (resp. two) straight lines parallel to the $x_1$-axis, see Figure~\ref{StandardExamples} right (resp. Figure~\ref{StandardExamples1} left).
In both cases, $R_D$ is the $\pi$-rotation around any of such line.

\begin{figure}
\begin{center}
\epsfysize=38mm 
\epsffile{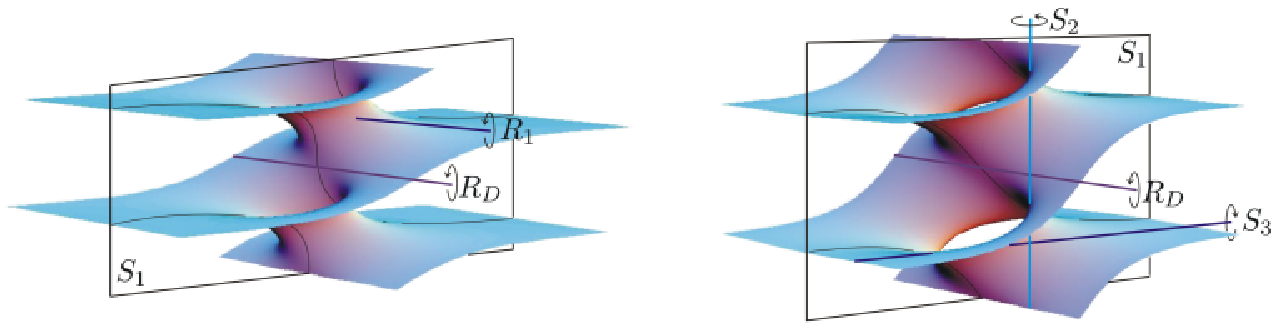}
\end{center}
\caption{Left: $M_{\t,0,\be}$ for $\t=\frac{\pi}{4}$ and $\be=\frac{3\pi}{8}$. 
Right: $M_{\t,0,\be}$ for $\t=\frac{\pi}{4}$ and $\be=\frac{\pi}{2}$.}
\label{StandardExamples1}
\end{figure}

\item In the case $\a=0$, $\be=\frac{\pi}{2}$, the puncture $A$ coincides with the lower right corner of $\widetilde{\cal R}$, and
Iso$(M_{\t,0,\frac{\pi}{2}})=$ Iso$(M_{\t,0,0})$.
The isometry $S_1$ represents in $\R^3$ a reflection symmetry across a plane orthogonal to the $x_1$-axis.
In this case, $S_3$ (resp. $S_2$) represents in $\R^3$ a $\pi$-rotation around one of the four (resp. two) straight lines parallel to the $x_2$-axis (resp. $x_3$-axis) contained on $M_{\t,0,\frac{\pi}{2}}$, see Figure~\ref{StandardExamples1} right.

\item If $0<\a<\frac{\pi}{2}$ and $\be=\frac{\pi}{2}$, then $S_3$ is an isometry of $(g,dh)$, since~$A$ moves from the lower right corner of $\widetilde{\cal R}$ to its upper right corner, as $\a$ varies from $0$ to $\frac{\pi}{2}$.
And Iso$(M_{\t,\a,\frac{\pi}{2}})$ is isomorphic to $(\Z /2\Z)^3$ with generators $S_3, {\cal D}, R_3= S_1\circ S_2$.
Now $S_3$ represents in $\R^3$ a $\pi$-rotation around any of the four straight lines parallel to the $x_2$-axis contained on $M_{\t,\a,\frac{\pi}{2}}$, and $R_3$ is the composition of a reflexion symmetry across a plane orthogonal to the $x_2$-axis with a translation by half a horizontal period, see Figure~\ref{MtaPi2}.

\begin{figure}
\begin{center}
\epsfysize=4cm 
\epsffile{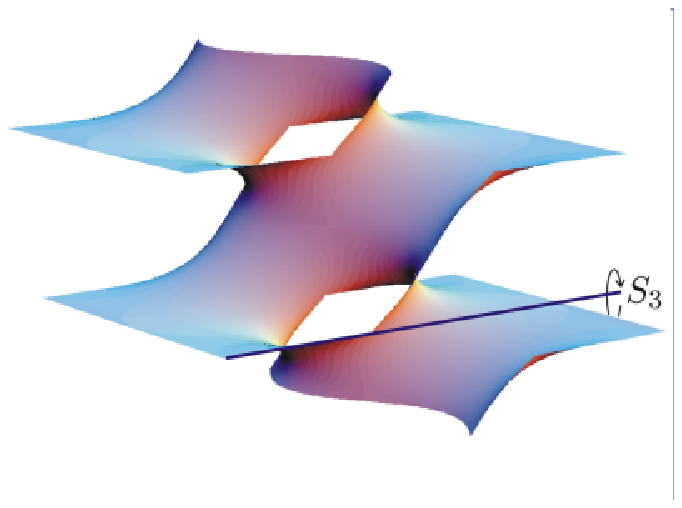}
\end{center}
\caption{The standard example $M_{\frac{\pi}{4},\frac{\pi}{4},\frac{\pi}{2}}$.}
\label{MtaPi2}
\end{figure}

\item Suppose now that $0<\a<\frac{\pi}{2}$ and $\be=0$.
The puncture $A$ moves horizontally to the right running along the upper boundary side of $\widetilde{\cal R}$.
Thus Iso$(M_{\t,\a,0})$ is isomorphic to $(\Z /2\Z)^3$, with generators $S_2$, ${\cal D}$, $R_2=S_1\circ S_3$.
As in the case of $M_{\t,0,0}$, $S_2$ represents in space a reflection symmetry across two planes orthogonal to the $x_2$-axis, and $R_2$ is a $\pi$-rotation around a line parallel to the $x_2$-axis that cuts $M_{\t,\a,0}$ orthogonally, see Figure~\ref{StandardExamples2} left.

\item If $\a =\frac{\pi}{2}$, then $M_{\t,\frac{\pi}{2},\be}$ is nothing but the rotated image of $M_{\t,\frac{\pi}{2},0}$ by angle $\be$ around the $x_3$-axis, hence we reduce the study to $\be=0$.
Now $A$ lies on the upper right corner of $\widetilde{\cal R}$, so all $S_1, S_2, S_3, R_D$ leave invariant the distribution of zeros and poles of the Gauss map of $M_{\t,\frac{\pi}{2},\be}$, and Iso$(M_{\t,\frac{\pi}{2},\be})=$Iso$(M_{\t,0,0})$.
$S_2$ represents a reflection symmetry across two planes orthogonal to the $x_2$-axis, $S_3$ (resp. $S_1$) represents in $\R^3$ a $\pi$-rotation around one of the four (resp. two) straight lines parallel to the $x_1$-axis (resp. $x_3$-axis) contained on $M_{\t,\frac{\pi}{2},0}$,
and $R_D$ corresponds to a reflection symmetry across two horizontal planes, see Figure~\ref{StandardExamples2} right.
\end{enumerate}

\mbox{}

Next let us show a uniqueness result we will use in the proof of Theorem~\ref{thmtopology} (Section~\ref{secKnotsimplyconnected}).
\begin{lemma}
\label{uniqueness}
With the notation above, $F_{\g_2}=(0,0,2\pi)$ if and only if $\a=\be=0$.
\end{lemma}
\begin{proof}
From (\ref{eq:periodhomologymuysim}) we know that $F_{\g_2}=(0,0,2\pi)$ when $\a=\be=0$.
Now suppose that $F_{\g_2}=(0,0,2\pi)$ and let us conclude that both $\a$ and $\be$ vanish.

If $\be=\frac{\pi}{2}$, $R_3$ is an isometry of $M_{\t,\a,\be}$, and we have $\g_2-(R_3)_*\g_2=\g_{A'}-\g_A$ and ${R_3^*\Phi=(\phi_1,-\phi_2,\phi_3)}$.
Moreover, we obtain from (\ref{eq:finales}) and (\ref{eq:periodendsnosim}) that $P_{\g_A}=P_{\g_{A'}}$ and
$-F_{\g_A}=F_{\g_{A'}}=(0,\pi a,0)$, with $a=\frac{\mu\sin\t}{\sqrt{1-\sin^2\t\cos^2\a}}>0$. Thus
\[
\int_{\g_2}\Phi=\int_{\g_2}(\phi_1,-\phi_2,\phi_3)+2 i(0,\pi a,0),
\]
and the second component of $F_{\g_2}$ equals $\pi a\neq 0$, which is not possible.
Hence it must be $\be\neq\frac{\pi}{2}$.
Since $M_{\t,\frac{\pi}{2},\be}$ differs from $M_{\t,\frac{\pi}{2},\frac{\pi}{2}}$ in a rotation about the $x_3$-axis,
it also holds $\a\neq\frac{\pi}{2}$.
This is, $\a,\be\in[0,\frac{\pi}{2})$, and we can choose
for every $\a,\be$ the same curve representative $\g_2$ as in the case $\a=\be=0$ (i.e. $\g_2=\{z\in\overline\C_1\ |\ |z|=1\}$).

Since $P_{\g_2}=(0,0,0)$, then $F_{\g_2}=(i\int_{\g_2} g dh,2\pi)\in\C\times\R$, and
\begin{equation}
\label{eq:flujovertical}
0=\int_{\g_2}g\, dh=
-2\mu \int_{-\pi}^\pi\frac{\cos\be\,\sin t+i(\sin\a\,\sin\be\,\sin t-\cos\a\,\cos t)}{|\sigma-\delta e^{-it}|^2\sqrt{\l^2+\l^{-2}+2\cos(2t)}}\, dt .
\end{equation}
Therefore, we deduce from $\Re(\int_{\g_2}g\, dh) = 0$ that
\[
\int_0^\pi\left(\frac{1}{|\sigma-\delta e^{-it}|^2}-\frac{1}{|\sigma-\delta e^{it}|^2}\right) \frac{\sin t}{\sqrt{\l^2+\l^{-2}+2\cos(2t)}}\, dt
\]
\[
=4\sin\be\int_0^\pi\frac{\sin^2 t}{|\sigma-\delta e^{-it}|^2 |\sigma-\delta e^{it}|^2\sqrt{\l^2+\l^{-2}+2\cos(2t)}}\, dt=0.
\]
The only possibility is then $\be=0$, and equation (\ref{eq:flujovertical}) reduces to
\[
2\mu i\cos\a \int_{-\pi}^\pi\frac{\,\cos t}{|\sigma-\delta e^{-it}|^2\sqrt{\l^2+\l^{-2}+2\cos(2t)}}\, dt=0,
\]
which is equivalent to
\[
\int_0^\pi\left(\frac{1}{|\sigma-\delta e^{-it}|^2}+\frac{1}{|\sigma-\delta e^{it}|^2}\right) \frac{\cos t}{\sqrt{\l^2+\l^{-2}+2\cos(2t)}}\, dt
\]
\[
=2\sin\a\int_0^\frac{\pi}{2}\frac{(\frac{1}{|\sigma-\delta e^{-it}|^2|\sigma+\delta e^{it}|^2}+\frac{1}{|\sigma-\delta e^{it}|^2|\sigma+\delta e^{-it}|^2})\cos t}{\sqrt{\l^2+\l^{-2}+2\cos(2t)}}\, dt=0,
\]
which only is satisfied for $\a=0$. Hence $\a=\be=0$, as we wanted to prove.
\end{proof}

\begin{figure}
\begin{center}
\epsfysize=5cm 
\epsffile{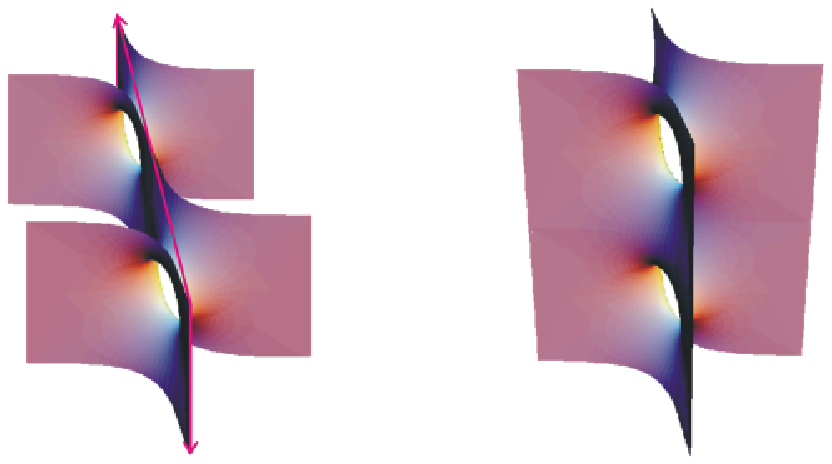}
\end{center}
\caption{Rotated image of $M_{\t,0,\be}$ (left) and two copies of the half of such a rotated $M_{\t,0,\be}$ (right), for $\t=\frac{\pi}{200}$ and $\be=\frac{9 \pi}{20}$.}
\label{LimitScherk1p}
\end{figure}

We finalize this subsection by listing all the possible degenerate limits of the standard examples $M_{\t,\a,\be}$, all of which can be directly computed by using the Weierstrass data.
\begin{itemize}
\item When $(\t,\be)\to(\t_0,\t_0)$, for some $\t_0\in(0,\frac{\pi}{2})$, $M_{\t,0,\be}$ converges smoothly to a Riemann minimal example.

\item When $\t \to 0^+$ and $(\a,\be)\to(0,0)$,
$M_{\t,\a,\be}$ converges smoothly to two catenoids with flux $(0,0,2\pi)$, see Figure~\ref{figureLimits} left.

\item When $\t \to 0^+$ and $(\a,\be)\to(\a_0,\be_0)\neq(0,0)$,
$M_{\t,\a,\be}$ converges smoothly to two copies of the singly periodic Scherk minimal surfaces with four ends, two of them horizontal, and with angle\footnote{We will call the {\it angle} of any singly or doubly periodic Scherk minimal surface to the angle between its nonparallel ends.} $\arccos(\cos \a_0 \cos \be_0)$, see Figure~\ref{LimitScherk1p}.

\item When $\t \to \frac{\pi}{2}^-$ and $(\a,\be)\to(0,\frac{\pi}{2})$, $M_{\t,\a,\be}$ converges smoothly (after blowing up) to two vertical helicoids spinning oppositely, see Figure~\ref{LimitHelicoid}.

\begin{figure}
\begin{center}
\epsfysize=2cm 
\epsffile{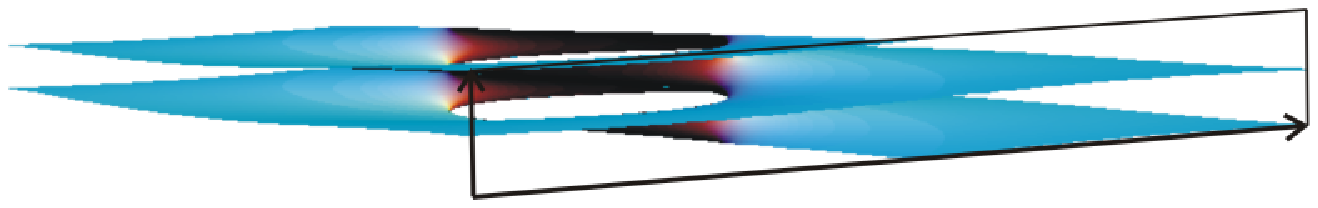}
\end{center}
\caption{$M_{\t,0,\be}$ for $\t=\frac{19 \pi}{40}$ and $\be=\frac{\pi}{2}$.}
\label{LimitHelicoid}
\end{figure}

\item When $\t \to \frac{\pi}{2}^-$ and $(\a,\be)\to(\a_0,\be_0)\neq(0,\frac{\pi}{2})$, $M_{\t,\a,\be}$ converges smoothly (after blowing up) to two copies of the doubly periodic Scherk minimal surfaces with four ends, two of them horizontal, and with angle $\arccos(\cos \a_0 \sin \be_0)$, see Figure~\ref{figureLimits} right.

\end{itemize}

\begin{remark}
\label{notaidentificaciones}
We have defined the $3$-parametric family $\{M_{\t,\a,\be}\ |\ (\t,\a,\be)\in {\cal I}_1\}$ of doubly periodic minimal surfaces, with ${\cal I}_1=\left\{(\t,\a,\be)\in\left(0,\frac{\pi}{2}\right) \times \left[0,\frac{\pi}{2}\right]^2\ |\ (\a,\be)\neq(0,\t)\right\}$, and in this range of parameters the end $A=A(\t,\a,\be)$ runs entirely along the closure of $\widetilde{\cal R}$ except its lower left corner $D'$.
$\widetilde{\cal R}$ can be identified conformally through the $z$-map with an octant of $\esf^2$.
We can easily extend the range of parameters so that $A$ runs entirely the sphere minus the branch values of the $z$-map, which can be achieved by varying $(\t,\a,\be)$ in ${\cal I}_2=\left\{(\t,\a,\be)\in(0,\frac{\pi}{2})\times[-\frac{\pi}{2},\frac{\pi}{2}]\times[-\pi,\pi]\ |\ (\a,\be)\neq(0,\pm\t),(0,\pm(\pi-\t))\right\}$.
We can define $M_{\t,\a,\be}$ for $(\t,\a,\be)\in{\cal I}_2$ similarly as for $(\t,\a,\be)\in{\cal I}_1$,
but it is straightforward to check that, up to translations, rotations and homotheties:
\begin{itemize}
\item $M_{\t,-\frac{\pi}{2},\be}$ coincides with $M_{\t,\frac{\pi}{2},\be}$, which does not depend on $\be$.
\item $M_{\t,-\a,0}$ is the reflected image of $M_{\t,\a,\be}$ with respect to a plane orthogonal to the $x_1$-axis.
\item $M_{\t,\a,\be\pm\pi}$ coincides with $M_{\t,\a,\be}$.
\item $M_{\t,0,-\be}$ is the reflected image of $M_{\t,\a,\be}$ with respect to a plane orthogonal to the $x_2$-axis.
\end{itemize}
\end{remark}

Therefore, we define the family of standard examples as ${\cal K}=\{M_{\t,\a,\be}\ |\ (\t,\a,\be)\in {\cal I}\}$, where
\begin{equation}\label{eq:I}
\textstyle{{\cal I}=\left\{(\t,\a,\be)\in(0,\frac{\pi}{2})\times(\frac{-\pi}{2},\frac{\pi}{2})\times[0,\pi)\ |\ (\a,\be)\neq(0,\t),(0,\pi-\t)\right\}\cup\{(\t,\frac{\pi}{2},0)\ |\ \t\in(0,\frac{\pi}{2})\}  .}
\end{equation}
We choose this space of parameters to avoid repeating surfaces twice, see Remark~\ref{notaidentificaciones}.

\begin{remark}
\label{remarkdim3}
By construction, the branch values of the Gauss map $N$ of $M_{\t,\a,\be}$ are contained in a spherical equator of $\esf^2$, so a consequence of Theorem~14 in~\cite{mro1} assures that the space of bounded Jacobi functions on $M$ coincides with the space of linear functions of $N$, $\{\langle N, V\rangle\ |\ V\in\R^3\}$ (in particular, such space is $3$-dimensional).
This condition is usually referred in literature as the {\it nondegeneracy} of $M_{\t,\a,\be}$, which can be interpreted by means of an Implicit Function Theorem argument to obtain that ${\cal K}$ is a $3$-dimensional real analytic manifold (Hauswirth and Traizet~\cite{HausTraizet1}).
\end{remark}

\subsection{The space of standard examples is self-conjugate}
\label{statement4}
In the previous subsection we have defined the family ${\cal K}=\{M_{\t,\a,\be}\ |\ (\t,\a,\be)\in {\cal I}\}$ of standard examples.
Given $M_{\t,\a,\be}\in{\cal K}$ with Weierstrass data $(g,dh)$, we let $M_{\t,\a,\be}^*$ denote the conjugate surface of $M_{\t,\a,\be}$, with Weierstrass data $(g,i dh)$.
Taking into account that the flux vector (resp. the period vector) of the conjugate surface along a given curve in the parameter domain equals the period vector (resp. the opposite of the flux vector) of the original surface along the same curve, we deduce from (\ref{eq:finales}), (\ref{eq:periodendsnosim}) and (\ref{eq:Perg2Mtab}) that $M_{\t,\a,\be}^*$ is a complete immersed torus invariant by the rank two lattice generated by the horizontal vector $P_{\g_A}^*=-F_{\g_A}$ and $P_{\g_2}^*=-F_{\g_2}$ (whose third coordinate is $-2\pi$) and which has four horizontal Scherk-type ends in the quotient.
Moreover, $M_{\t,\a,\be}^*$ is embedded thanks to the Maximum Principle, since the heights of its ends depend continuously on $(\a,\be)$, and $M_{\t,0,0}^*$ is embedded (it is constructed from congruent blocks being Jenkins-Serrin graphs).

Note that, by (\ref{eq:Perg2Mtab}), the period vector of $M_{\t,\a,\be}^*$ along $\g_2^*=\g_1+\g_A$ vanishes, and that the third component of the flux of $M_{\t,\a,\be}^*$ along $\g_2^*$ equals $f_1(\t)$ given by~(\ref{eq:f1}).
The next lemma ensure that, after scaling and rotating the surfaces around the $x_3$-axis, the families ${\cal K}$ and ${\cal K}^*=\{M_{\t,\a,\be}^*\ |\ (\t,\a,\be)\in {\cal I}\}$ coincide, which finishes the proof of Theorem~\ref{thmK}.

\begin{lemma}
\label{lemaautoconj}
Given $(\t,\a,\be)\in {\cal I}_1$, the surface $M_{\frac{\pi}{2}-\t,\a,\be+\frac{\pi}{2}}$ coincides with $M_{\t,\a,\be}^*$ up to normalizations.
\end{lemma}
\begin{proof}
It is easy to see that $\Sigma_{\frac{\pi}{2}-\t}=\{(\widetilde{z},\widetilde{w})\ |\ \widetilde{w}^2=(\widetilde{z}^2-1)^2+4\widetilde{z}^2\sec^2\t\}$.
Since the M\"{o}bius transformation $\varphi(z)=\frac{1-iz}{z-i}$ takes the set of branch points of the $z$-projection of $\Sigma_\t$ bijectively to the set of branch points
of the $\widetilde{z}$-projection of $\Sigma_{\frac{\pi}{2}-\t}$, it follows that $\Theta(z,w)=(\varphi(z),\widetilde{w}(\varphi(z)))$ is a biholomorphism between $\Sigma_\t$ and $\Sigma_{\frac{\pi}{2}-\t}$.
On the other hand, it is straightforward to check that $g_{\t,\a,\be}=g_{\frac{\pi}{2}-\t,\a,\be-\frac{\pi}{2}}
\circ \Theta $, where the subindex means the parameters of the standard example $M_{\t,\a,\be}$ for which the corresponding $g_{\t,\a,\be}$ is the Gauss map. Denoting by $dh_{\t}$ its height differential (recall that it only depends on $\t $)
a direct computation gives $\Theta^*dh_{\frac{\pi}{2}-\t}
=-\frac{\mu(\frac{\pi}{2}-\t)}{\mu(\t)\tan \t}\, i dh_{\t}
=-\frac{K(\sin^2\t)}{K(\cos^2\t)}\, i dh_{\t}$.
Hence $M_{\frac{\pi}{2}-\t,\a,\be-\frac{\pi}{2}}=M_{\t,\a,\be}^*$ up to normalizations.
Since $M_{\frac{\pi}{2}-\t,\a,\be-\frac{\pi}{2}}=M_{\frac{\pi}{2}-\t,\a,\be+\frac{\pi}{2}}$ by Remark~\ref{notaidentificaciones}, the lemma is proved.
\end{proof}

\section{The classifying map}
\label{secclassifyingmap}
The surfaces in ${\cal K}$ can be naturally seen inside the 4-dimensional complex manifold ${\cal W}$ consisting roughly  of all admissible Weierstrass data in the setting of Theorem~\ref{thmuniqueness}.
\begin{definition}
\label{defmarked}
We denote by ${\cal W}$ the space of tuples $(\M,g;p_1,p_2,q_1,q_2,[\g])$, where $g$ is a degree two meromorphic map defined on a torus $\M$ which is unbranched at its zeros $p_1,p_2$ and at its poles $q_1,q_2$, and $[\g]$ is a homology class in $\M-\{p_1,p_2,q_1,q_2\}$, which is not trivial in $H_1(\M,\Z)$.
\end{definition}
See~\cite{PeRoTra1} for a detailed description of ${\cal W}$.
We will shorten the elements in ${\cal W}$ simply by $g$, and call them {\it marked meromorphic maps}.
Each $g=(\M,g;p_1,p_2,q_1,q_2,[\g])\in{\cal W}$ determines a unique holomorphic differential $\phi=\phi(g)$ on $\M$ such that
\begin{equation}
\label{eq:defphi} \int_\g \phi =2\pi i,
\end{equation}
since the complex space of holomorphic differentials on $\M$ has dimension one.
Thus each $g\in{\cal W}$ can be seen as the Weierstrass data $(g,\phi)$, defined on $g^{-1}(\C^*)$, of a potential surface in the setting of Theorem~\ref{thmuniqueness}.
Equation~(\ref{eq:defphi}) means that the period vector of $(g,\phi)$ along $\g$ is horizontal and its flux along $\g$ has
third coordinate $2\pi$.

\begin{definition}We will say that $g\in{\cal W}$ {\it closes periods} when the next equations hold
\begin{equation}
\label{eq:gclosesperiods}
\int_\g \frac{\phi}{g}=\overline{\int_\g g\, \phi} \quad\mbox{and }\quad
\mbox{Res}_{p_1}\frac{\phi}{g}=-\mbox{Res}_{q_1}(g\, \phi)= a,\ \mbox{for certain } a\in\R^+.
\end{equation}
\end{definition}

Note that the first equation in (\ref{eq:gclosesperiods}), together with (\ref{eq:defphi}),
says that $P_\g=(0,0,0)$ and $F_\g=(i\int_\g g\, \phi, 2\pi)\in \C\times\R$.
The next lemma justifies the above definition of closing periods.

\begin{lemma}[\cite{PeRoTra1}]
\label{lemagclosesperiods}
If $g\in{\cal W}$ closes periods, then $(g,\phi)$ is the Weierstrass pair of a properly immersed minimal surface $M\subset \T\times\R$, for a certain flat torus $\T$, with total curvature $8\pi$ and four horizontal Scherk-type ends.
Furthermore, the fluxes at the ends $p_j,q_j$ are equal to $(-1)^{j+1}(\pi a,0,0)$ for the positive real number $a$ appearing in~(\ref{eq:gclosesperiods}), $j=1,2$.
\end{lemma}

Next we describe how to see each standard example $M_{\t,\a,\be}$ as an element of ${\cal W}$ which closes periods.
In a first step we rotate $M_{\t,\a,\be}$ about the $x_3$-axis so that the period $P_{\g_A}$ at its end $A$ (we follow the notation in Section~\ref{secstandardexamples}) is $(0,\pi a,0)$ for certain $a>0$ .
Now we associate to $M_{\t,\a,\be}$ the marked meromorphic map
\[
(\Sigma_\t,g;A'''={\cal F}(A),A'={\cal E}(A),A,A''={\cal D}(A),[\g_2]),
\]
where everything has been already defined in Subsection~\ref{subsecMtab} except the homology class~$[\g_2]$.
Recall that the ends $A,A',A'',A'''$ depend continuously on $\a,\be$ and that we described explicitly the loop $\g_2$ for $\a=\be=0$. For the remaining values of $\a,\be$, we take a embedded closed curve $\g_2\subset\Sigma_\t-\{A,A',A'',A'''\}$ depending continuously on $\a,\be$ so that $[\g_2]$ remains  constant in $H_1(\Sigma_\t,\Z)$.

\subsection{The ligature map}
We call {\it ligature map} to the holomorphic map $L:{\cal W}\to \C^4$ defined as follows
\[
L(g)=\left(\mbox{Res}_{p_1}\frac{\phi}{g}, \mbox{Res}_{q_1}(g\, \phi), \int_\g\frac{\phi}{g}, \int_\g g\, \phi\right),
\]
which clearly distinguishes when a marked meromorphic map closes periods:
\begin{quote}
{\it A marked meromorphic map $g\in{\cal W}$ closes periods if and only if there exist $a\in\R^+$ and $b\in\C$ such that $L(g)=(a,-a,\overline{b},b)$.}
\end{quote}
Since the residues of a meromorphic differential on a compact Riemann surface add up to zero, if the second equation in (\ref{eq:gclosesperiods}) holds, then $\mbox{Res}_{p_2}\frac{\phi}{g}=-\mbox{Res}_{q_2}(g\, \phi)= -a$.

Let ${\cal S}=\{S_\rho\ |\ \rho\in(0,\pi)\}$ be the 1-dimensional moduli space of singly periodic Scherk minimal surfaces with two horizontal ends, vertical part of the flux at its two nonhorizontal ends equals to $2\pi$ and period vector in the direction of the $x_2$-axis.
For each $\rho\in(0,\pi)$, let $S_\rho\in{\cal S}$ the singly periodic Scherk surface of angle $\rho$. The limiting normal vectors of $S_\rho$ at its nonhorizontal ends project stereographically to $\tan\frac{\rho}{2},-\cot\frac{\rho}{2}$.
Recall we can obtain two copies of $S_\rho$ by taking limits from standard examples, see Subsection~\ref{subsecMtab}.
We identify $S_\rho$ with the list $(\M_\rho, g;0_1,0_2,\infty_1,\infty_2,[\g_\rho])$, where:
\begin{itemize}
\item $\M_\rho$ is a Riemann surface with nodes constructed by gluing two copies $\overline\C_1,\overline\C_2$ of $\overline\C$ with nodes $\tan\frac{\rho}{2},-\cot\frac{\rho}{2}$.
\item $g:\M_\rho\to\overline{\C}$ is the  map which associates to each point in $\M_\rho$ its complex value as a point in $\overline\C_j$, $j=1,2$
(in particular, the degree of $g$ equals two).
\item $0_j,\infty_j$ are respectively the zero and infinity in $\overline\C_j$, $j=1,2$.
\item $\g_\rho\subset\C_1$ is an embedded closed curve in the homology class $[\Gamma_1]+[\Gamma_2]$, where $\Gamma_1$
(resp.~$\Gamma_2$) is a small loop in $\overline\C_1$ around $0_1$ (resp. $\tan\frac{\rho}{2}$) with the positive orientation.
\end{itemize}

\begin{figure}
\begin{center}
\epsfysize=44mm 
\epsffile{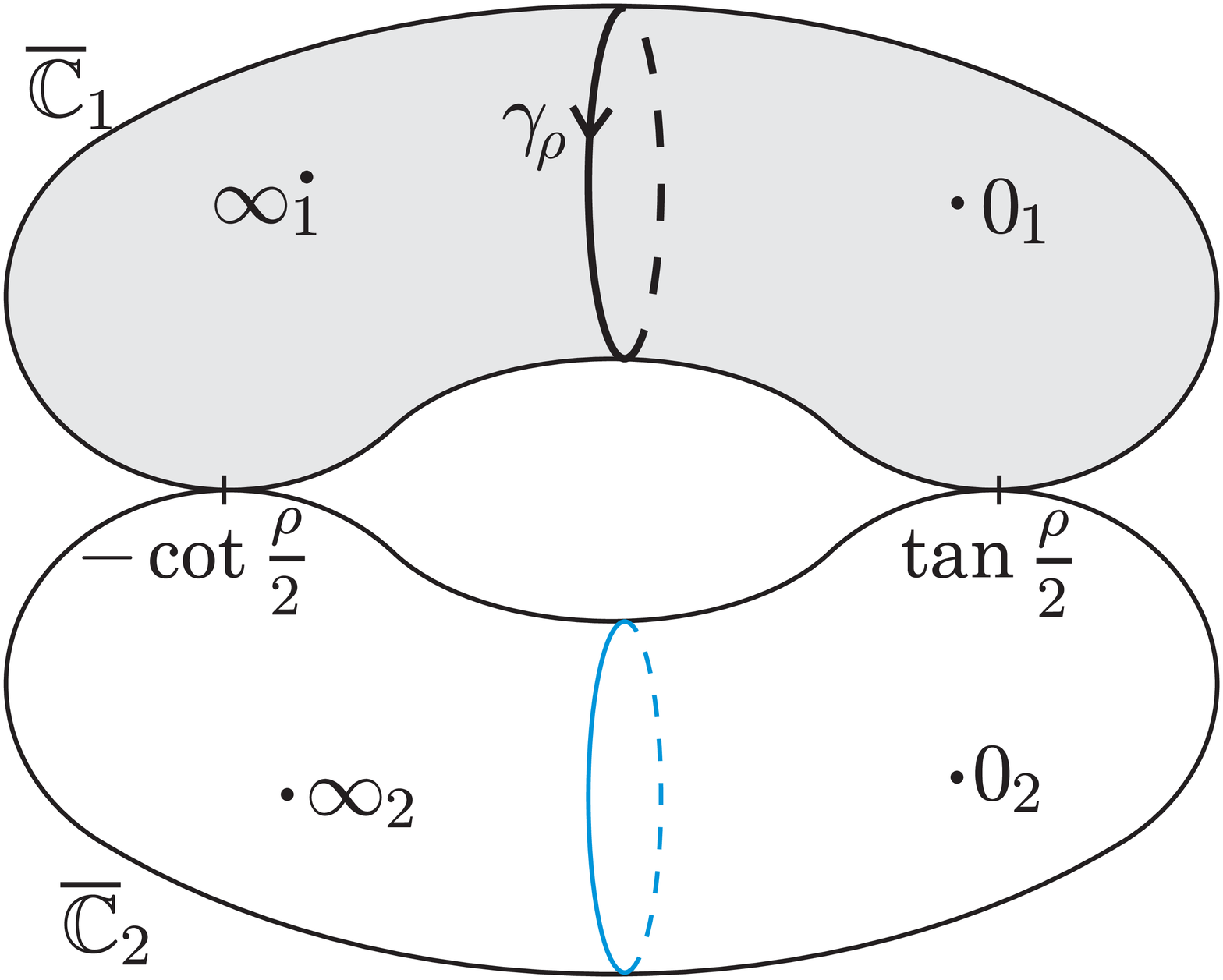}
\end{center}
\caption{The Riemann surface $\M_\rho$ with nodes $\tan\frac{\rho}{2},-\cot\frac{\rho}{2}$, 
and the embedded closed curve $\g_\rho\subset\overline\C_1$.}
\label{morcillas}
\end{figure}

With the above identification, we can see ${\cal S}$ in $\partial{\cal K}\subset\partial{\cal W}$.
From Lemma~6 and Theorem~2 of~\cite{PeRoTra1}, one deduces that $\widetilde{\cal W}={\cal W}\cup{\cal S}$ is a $4$-dimensional complex manifold where $L$ extends holomorphically, and this extended ligature map is a biholomorphism in a small neighborhood of ${\cal S}$ in $\widetilde{\cal W}$.
A straightforward computation gives~us
\begin{equation}
\label{eq:LS}
\textstyle{
L(S_\rho)=
\left(2\csc\rho,-2\csc\rho,-2\pi i\tan\frac{\rho}{2},2\pi i\tan\frac{\rho}{2}\right) .}
\end{equation}

\subsection{The classifying map}
In this subsection we will study the topology of the space ${\cal K}$, and the key ingredient for this study will be the classifying map $C$ that associates roughly to each marked standard surface its period at the ends and the horizontal component of its flux along a nontrivial homology class with zero period vector.
\begin{definition}
{\rm We denote $\widetilde{\cal K}={\cal K}\cup{\cal S}$ and define the {\it classifying map} $C:\widetilde{\cal K}\to \R^+\times \C$ by $C(M)=(a,b)$ so that $\mbox{\rm Res}_{p_1}\frac{dh}{g}=a$ and $\int_\g g\, dh=b$ (hence $F_{\g}=(i b,2\pi)\in\C\times\R$), provided that $M=(\M,g;p_1,p_2,q_1,q_2,[\g])$ and $dh$ is the height differential of $M$.}
\end{definition}

Let $M=(\M,g;p_1,p_2,q_1,q_2,[\g])$ be a marked surface in $\widetilde{\cal K}$, and $C(M)=(a,b)$.
Denote by $\g_X$ a small loop around $X\in\M$, oriented positively.
If $\widetilde M=(\M,g;p_1,p_2,q_1,q_2,[\widetilde\g])$
with $[\widetilde\g]=[\g]+n([\g_{p_1}]+[\g_{q_1}])+m([\g_{p_2}]+[\g_{q_2}])$, then
$C(\widetilde M)=(a,b+2\pi a(n-m))$, see Lemma~\ref{lemagclosesperiods}.
Since we want to avoid associating more than one different image to the same geometrical surface, it is necessary to restrict
\[
C:\widetilde{\cal K}\to \Lambda=\{(a,b)\in\R^+\times\C\ |\ 0\leq\Re(b)<2\pi a\}\equiv\R^+\times\esf^1\times\R .
\]

Recall that ${\cal K}$ is a 3-dimensional real analytic manifold, see Remark~\ref{remarkdim3}.
It is clear from the definition above that $C|_{\cal K}$ is smooth. Also note that $C$ is essentially $L|_{\widetilde{\cal K}}$.
Since $L$ extends as a biholomorphism in a small neighborhood of ${\cal S}$ in $\widetilde{\cal W}$, it follows that $\widetilde{\cal K}$ can be endowed with a structure of a 3-dimensional real analytic manifold and $C:\widetilde{\cal K}\to \Lambda$ is also smooth.
Remark that $C|_{\cal K}$ is not proper, since ${\cal S}$ is contained in the boundary of ${\cal K}$ in $\widetilde{\cal W}$ but $C({\cal S})\subset\Lambda$.
In detail, from (\ref{eq:LS}) we have
\begin{equation}
\label{eqC(S)}
C(S_\rho)=\left(2\csc\rho,2\pi i\tan\frac{\rho}{2}\right) .
\end{equation}

\begin{proposition}
\label{proposproper}
The classifying map $C:\widetilde{\cal K}\to \Lambda$ is proper.
\end{proposition}
\begin{proof}
Take a sequence $\{M_n\}_n\subset\widetilde{\cal K}$ so that $\{C(M_n)=(a_n,b_n)\}_n$ converges to some point $(a,b)\in\Lambda$,
and let us prove that a subsequence of $\{M_n\}_n$ converges to a surface in $\widetilde{\cal K}$.

First suppose that, after passing to a subsequence, $M_n\in{\cal K}$ for every $n$, and let $(\t_n,\a_n,\be_n)\in{\cal I}$,
see (\ref{eq:I}), be the angles which determine the spherical configuration of $M_n=M_{\t_n,\a_n,\be_n}$.
Extracting a subsequence, we can assume that $(\t_n,\a_n,\be_n)\to(\t_\infty,\a_\infty,\be_\infty)\in[0,\frac{\pi}{2}]\times[-\frac{\pi}{2},\frac{\pi}{2}]\times[0,\pi]$.
By equation (\ref{eq:periodendsnosim}) we deduce that $a(M_n)$ equals the modulus of $\mu(\t_n) \sin\t_n\, E(\t_n,\a_n,\be_n)\in\C$, and so:
\begin{itemize}
\item If $\t_\infty=\frac{\pi}{2}$, then $a(M_n)\to 0$.
These limits correspond, after blowing up, to two vertical helicoids when $\a_\infty=0$ and $\be=\frac{\pi}{2}$, or to two copies of a doubly periodic Scherk minimal surface otherwise, see Subsection~\ref{subsecMtab}.
\item If $\t_\infty\neq\frac{\pi}{2}$, but $\a_\infty=0$ and $\be_\infty\in\{\t_\infty, \pi-\t_\infty\}$, then $a(M_n)\to\infty$.
These limits correspond to the vertical catenoid when $\t_\infty=0$, or to a Riemann minimal example otherwise.
\end{itemize}
Therefore, the only possibilities are:
\begin{itemize}
\item $\t_\infty=0$ and $(\a_\infty,\be_\infty)\not\in\{(0,0),(0,\pi)\}$, and then  $\{M_n\}_n$ converges to two copies of a singly periodic Scherk minimal surface.
\item $(\t_\infty,\a_\infty,\be_\infty)\in{\cal I}$, so $\{M_n\}_n$ converges to the standard example $M_{\t_\infty,\a_\infty,\be_\infty}\in{\cal K}$.
\item $\a_\infty=\pm\frac{\pi}{2}$, hence $M_n\to M_{\t_\infty,\pm\frac{\pi}{2},\be_\infty}=M_{\t_\infty,\frac{\pi}{2},0}$ (see Remark~\ref{notaidentificaciones}).
\item $\be_\infty=\pi$, hence $M_n\to M_{\t_\infty,\a_\infty,\pi}=M_{\t_\infty,\a_\infty,0}$.
\end{itemize}
Hence $\{M_n\}_n$ admits a subsequence converging in $\widetilde{\cal K}$
(this can be also obtained by arguing as in the proof of Theorem~5 in~\cite{PeRoTra1}).

Thus it suffices to prove that $C|_{\cal S}$ is proper, but this is clear by~(\ref{eqC(S)}).
This fact completes the proof of Proposition~\ref{proposproper}.
\end{proof}

\section{The classifying map is a local diffeomorphism}
\label{secDiff}
\begin{proposition}
\label{ThDiff}
The classifying map $C:\widetilde{\cal K}\to \Lambda$ is a local diffeomorphism.
\end{proposition}
\begin{proof}
The relationship between $C$ and $L|_{\widetilde{\cal K}}$ allows us to assure that $C$ is a diffeomorphism in a small neighborhood of ${\cal S}$ in $\widetilde{\cal K}$.
Thus it only remains to demonstrate that $C|_{\cal K}$ is a local diffeomorphism.
Consider a standard example $M\in{\cal K}$ and denote by $\widetilde{M}$ its lifting to $\R^3$.
It suffices to check that if $u:\widetilde{M}\to\R$ is a Jacobi function that lies in the kernel of $d C_M$, then $u=0$.

We can write $u=\langle\left.\frac{d}{dt}\right|_0\widetilde{M}_t,N\rangle$ for certain variation $\{\widetilde M_t\}\subset{\cal K}$ of $\widetilde{M}_{t=0}=\widetilde{M}$.
Let ${\cal P}_t$ be the period lattice of $\widetilde{M}_t\subset\R^3$, $M_t=\widetilde{M}_t/{\cal P}_t$ and $(a_t,b_t)=C(M_t)\in\R^+\times\C$. Since $u\in\ker(d C_M)$, then
\begin{equation}
\label{eqder}
\textstyle{\left.\frac{d}{dt}\right|_{t=0}a_t=0\qquad \mbox{and }\qquad \left.\frac{d}{dt}\right|_{t=0}b_t=0.}
\end{equation}

By Lemma~\ref{lemaautoconj} and after normalizations, the conjugation map $*:{\cal K}\to{\cal K}$, which associates to each standard example its conjugate surface, is a well-defined map.
Furthermore, $*$ is clearly differentiable since it is the restriction to ${\cal K}$ of the map $(g,\phi)\to(g,i \phi)$ on the space of allowed Weierstrass data.
Since clearly $*\circ *=$identity, then we deduce that $*$ is a diffeomorphism.
Hence it suffices to prove that the tangent vector $v$ defined as the image of $u$ by the differential of $*$, vanishes identically.
Notice that $v=\langle\left.\frac{d}{dt}\right|_{t=0}\widetilde{M}^*_t,N\rangle$, being $\widetilde{M}^*_t$ the conjugate surface of $\widetilde{M}_t$.
In particular, $v$ is a Jacobi function on $\widetilde{M}^*$, which is moreover bounded since all the $\widetilde{M}^*_t$ have horizontal ends.

First suppose that $v$ is a bounded Jacobi function on the quotient $M^*$ of $\widetilde{M}^*$ by its period lattice.
By Remark~\ref{remarkdim3} we know that $v$ is of the kind $v=\langle N,V\rangle$, for some $V\in\R^3$.
Theorem~3 in~\cite{mro1} assures that there exists a unique element $X_v$ of the space of complete branched minimal immersions into $\R^3$ (including the constant maps) with finite total curvature and planar ends whose extended Gauss map is $N$, such that $\langle X_v,N\rangle=v$.
Thus $X_v$ is constantly $V$, and $v$ corresponds to a translation of $\widetilde{M}^*$ in $\R^3$.
As we are considering the surfaces in ${\cal K}$ up to translations, it holds $v=0$.
Therefore it only remains to prove that $v$ descends to the quotient $M^*$.

Recall that the flux of $\widetilde{M}_t$ at its ends is (up to sign) $H_t=(\pi a_t,0,0)$ and its flux along the homology class in the last component (viewed as a marked standard example) is $T_t=(i b_t,2\pi)\in\C\times\R\equiv\R^3$.
Therefore the period lattice of $\widetilde{M}^*_t$ (before normalization) is generated by $H_t$ and $T_t$.
We parameterize $\widetilde{M}^*_t$ by $\psi^*_t:\widetilde M\to\widetilde{M}^*_t$,
and denote by $S_{1,t},S_{2,t}:\widetilde{M}\to\widetilde{M}$ the diffeomorphisms induced by $H_t,T_t$, i.e. those satisfying
\begin{equation}
\label{eqperiodos}
\psi^*_t\circ S_{1,t}=\psi^*_t+H_t \qquad \mbox{and }\qquad
\psi^*_t\circ S_{2,t}=\psi^*_t+T_t .
\end{equation}
By (\ref{eqder}), $\left.\frac{d}{dt}\right|_{t=0}H_t=\left.\frac{d}{dt}\right|_{t=0}T_t=\vec{0}$.
Therefore,
\[
\begin{array}{rcl}
v\circ S_{1,0}& = & \langle\left.\frac{d}{dt}\right|_{t=0}\psi^*_t,N\rangle\circ S_{1,0}
=\langle\left.\frac{d}{dt}\right|_{t=0}(\psi^*_t\circ S_{1,t}),N\rangle\\
 & \stackrel{(\ref{eqperiodos})}{=}&
\langle\left.\frac{d}{dt}\right|_{t=0}\psi^*_t,N\rangle+
\langle\left.\frac{d}{dt}\right|_{t=0}H_t,N\rangle =v .\\
\end{array}
\]
Analogously $v\circ S_{2,0}=v$. Thus $v$ descends to the quotient and Proposition~\ref{ThDiff} is proved.
\end{proof}

\section{The topology of ${\cal K}$ (proof of Theorem~\ref{thmtopology})}
\label{secKnotsimplyconnected}
By Propositions~\ref{proposproper} and~\ref{ThDiff}, $C:\widetilde{\cal K}\to \Lambda$ is a proper local diffeomorphism, and so a finite sheeted covering map. We deduce from Lemma~\ref{uniqueness}, Remark~\ref{notaidentificaciones} and equation (\ref{eqC(S)}) that the only surfaces $M\in\widetilde{\cal K}$ with $C(M)=(a,0)$, for some $a>0$, are the standard examples $M_{\t,0,0}$.
Since $a(M_{\t,0,0})=\mu(\t)$ is a strictly decreasing function in $\t$, the number of sheets of the covering map $C$ is one, and so it is a diffeomorphism.
Moreover, the set $C({\cal S})$ consists of the proper arc
$\rho\in(0,\pi)\mapsto \left(2\csc\rho,2\pi i\tan\frac{\rho}{2}\right)$,
from where it follows that ${\cal K}$ is diffeomorphic to the complement in $\Lambda$ of such arc, which in turn is diffeomorphic to $\R\times(\R^2-\{(\pm 1,0)\})$.
This proves Theorem~\ref{thmtopology}.

\begin{remark}
Since $C:{\cal K}\to C({\cal K})\equiv \R\times(\R^2-\{(\pm 1,0)\})$ is a diffeomorphism, any standard example is completely determined by its image through $C$. This justifies the words {\em classifying map} for~$C$.
\end{remark}

Recall that we can identify some standard examples in $\cal K$ by symmetries (see Remark~\ref{notaidentificaciones}).
Since there are standard examples invariant by such symmetries, we deduce that the quotient of $\cal K$ by these symmetries has structure of $3$-dimensional orbifold. Notice that this quotient space is the moduli space of doubly periodic minimal surfaces with parallel ends and genus one in the quotient.

\center{M. Magdalena Rodr\'\i guez at magdarp@ugr.es}

\bibliographystyle{plain}

\begin{thebibliography}{10}

\bibitem{HausTraizet1}
L.~Hauswirth and M.~Traizet.
\newblock The space of embedded doubly-periodic minimal surfaces.
\newblock {\em Indiana Univ. Math. J.}, 51(5):1041--1079, 2002.
\newblock MR1947868.

\bibitem{hm10}
D.~Hoffman and W.~H. Meeks~III.
\newblock The strong halfspace theorem for minimal surfaces.
\newblock {\em Invent. Math.}, 101:373--377, 1990.
\newblock MR1062966, Zbl 722.53054.

\bibitem{ka4}
H.~Karcher.
\newblock Embedded minimal surfaces derived from {S}cherk's examples.
\newblock {\em Manuscripta Math.}, 62:83--114, 1988.
\newblock MR0958255, Zbl 658.53006.

\bibitem{ka6}
H.~Karcher.
\newblock Construction of minimal surfaces.
\newblock {\em Surveys in Geometry}, pages 1--96, 1989.
\newblock University of Tokyo, 1989, and Lecture Notes No. 12, SFB256, Bonn,
  1989.

\bibitem{ml1}
H.~Lazard-Holly and W.~H. Meeks~III.
\newblock Classification des surfaces minimales de genre z\'ero proprement
  plong\'ees dans $\rth / {{\mathbb Z} ^2 } $.
\newblock {\em Comptes--Rendus de l'Acad{\'e}mie des Sciences de Paris}, pages
  753--754, 1997.

\bibitem{mpr1}
W.~H. Meeks~III, J.~P\'{e}rez, and A.~Ros.
\newblock Uniqueness of the {R}iemann minimal examples.
\newblock {\em Invent. Math.}, 131:107--132, 1998.
\newblock MR1626477, Zbl 916.53004.

\bibitem{mr3}
W.~H. Meeks~III and H.~Rosenberg.
\newblock The global theory of doubly periodic minimal surfaces.
\newblock {\em Invent. Math.}, 97:351--379, 1989.
\newblock MR1001845, Zbl 676.53068.

\bibitem{mro1}
S.~Montiel and A.~Ros.
\newblock Schr\"{o}dinger operators associated to a holomorphic map.
\newblock In {\em Global Differential Geometry and Global Analysis (Berlin,
  1990)}, volume 1481 of {\em Lecture Notes in Mathematics}, pages 147--174.
  Springer-Verlag, 1991.
\newblock MR1178529, Zbl 744.58007.

\bibitem{PeRoTra1}
J.~P\'{e}rez, M.~M. Rodr\'{\i}guez, and M.~Traizet.
\newblock The classification of doubly periodic minimal tori with parallel
  ends.
\newblock {\em J. of Differential Geometry}, 69(3):523--577, 2005.
\newblock MR2170278, Zbl pre05004289.

\bibitem{sche1}
H.~F. Scherk.
\newblock Bemerkungen \"{u}ber die kleinste {F}l\"{a}che innerhalb gegebener
  {G}renzen.
\newblock {\em J. R. Angew. Math.}, 13:185--208, 1835.

\end{thebibliography}

\end{document}